\documentclass[english]{ourlematema}
\usepackage{amsmath, amsthm, amssymb, hyperref, color}
\usepackage{graphicx}
\usepackage{caption}
\usepackage{verbatim}
\usepackage{chemfig, chemnum}
\usepackage{tikz}
\usepackage[linesnumbered,lined,commentsnumbered]{algorithm2e}
\usepackage{caption}
\usepackage{stmaryrd}
\usepackage{subcaption}
\usepackage{float}
\usepackage{mathtools}
\usepackage{braket}
\usepackage{anyfontsize}
\usepackage{todonotes}
\usepackage{tikz-cd}
\usepackage{cleveref}
\usepackage{csquotes}
\usepackage{enumitem}
\setenumerate{label=\textup{(\roman*)}}

\newtheorem{examplex}[thm]{Example}
\renewenvironment{exa}
{\pushQED{\qed}\examplex}
{\popQED\endexamplex}
\definecolor{ourgreen}{rgb}{0, 0.7, 0}
\hypersetup{
    colorlinks=true,
    linkcolor=blue,
    filecolor= green,      
    urlcolor=blue,
    citecolor=ourgreen,
}

 \title{What is Positive Geometry?}
 
  \author{Kristian Ranestad}
  \address{%
  University of Oslo, Oslo \\
\email{ranestad@math.uio.no}
}

  \author{Bernd Sturmfels}
  \address{%
  MPI for Mathematics in the Sciences, Leipzig \\
\email{bernd@mis.mpg.de}
}

\author{Simon Telen}
\address{%
  MPI for Mathematics in the Sciences, Leipzig \\
\email{simon.telen@mis.mpg.de}
}

\begin{document}

\maketitle

\begin{abstract}
    \noindent This article serves as an introduction to the special volume on Positive Geometry in the journal Le Matematiche. We attempt to answer the question in the title by describing the origins and objects of positive geometry at this early stage of its development. We discuss the problems addressed in the volume and report on the progress. We also list some open challenges. 
\end{abstract}

\vspace{-0.6cm}

\section{Introduction}

The volume \emph{Positive Geometry} collects a series of articles which report on collaborations that were initiated during the school \emph{Combinatorial Algebraic Geometry from Physics}, held at the Max Planck Institute for Mathematics in the Sciences, Leipzig, in May 2024. Several projects were inspired by the lectures of Michael Borinsky and Thomas Lam at that workshop. The topics of the articles align with the current goals of positive geometry, hence the volume title. Though it might be too early to spell out a definition of the field, the present article explains the editors' view on what is currently meant by positive geometry.

Positive geometry is an interdisciplinary subject inspired by fundamentally new ideas in particle physics and cosmology. Its origins lie in the observation that certain physical observables, called \emph{scattering amplitudes} and \emph{cosmological correlators}, can be recovered from the geometry of elegant mathematical objects, such as \emph{amplituhedra} \cite{arkani2014amplituhedron} and \emph{cosmological polytopes} \cite{arkani2017cosmological}. This has motivated physicists to search for such geometric structures in a range of quantum field theories. In positive geometry, we put this on a solid mathematical footing. 

A first step in this direction was the formal definition of
\emph{positive geometries} in \cite{arkani2017positive}. Building
\emph{Positive Geometry} as a subject is the next step.
A positive geometry is a tuple consisting of a complex algebraic variety $X$, a semi-algebraic subset $X_{\geq 0}$ of its real points, and a meromorphic top-form $\Omega(X_{\geq 0})$. This tuple satisfies recursive axioms prescribed by physics. These are listed in \cite[Section 2.1]{arkani2017positive}, \cite[Definition 1]{lam2024invitation} and \cite[Definition 1.20]{Lam}. 
Here, rather than spelling out the axioms, we present examples highlighting the topics featured in 
this~volume. 

\begin{exa} \label{ex:introduction}
    Let $X_{\geq 0} = P \subset \mathbb{R}^n$ be a convex polytope of dimension $n$. For an interior point $x \in {\rm int}(P)$, let $(P-x)^\circ$ be the polar dual of the translated polytope $P-x$. The volume of $(P-x)^\circ$ depends on $x$. It defines a meromorphic~form 
    \[ \Omega(P) \, = \, {\rm vol}(P-x)^\circ \, {\rm d}x\]
    on $X = \mathbb{P}^n$, called the canonical form of $P$ in \cite{arkani2017positive}. The poles of $\Omega(P)$ are simple, and they coincide with the arrangement of hyperplanes  defined by the facets of $P$. The zeros of $\Omega(P)$ lie along the adjoint hypersurface of $P$. If $P$ is a cosmological polytope, then the cosmological correlator corresponding to $P$ is $\int_{\Gamma} \Omega(P)$ for a suitable integration cycle $\Gamma$. In that context, the facet hyperplanes of $P$ depend on parameters, here denoted by $z$. The function $z \mapsto \int_{\Gamma} \Omega(P(z))$ satisfies linear differential equations, best encoded via holonomic $D$-modules. 

     In particle physics, one replaces $\mathbb{P}^n$ with a complex Grassmannian variety $X = {\rm Gr}_{\mathbb{C}}(k,n)$, and the polytope $P$ by a semi-algebraic subset ${\cal A} = X_{\geq 0}$ of its real points, called the (tree) amplituhedron. Scattering amplitudes in $N=4$ Super Yang--Mills theories are integrals of the canonical differential form $\Omega({\cal A})$. In this setting, the facet hyperplane arrangement is replaced with a divisor in ${\rm Gr}_{\mathbb{C}}(k,n)$, obtained as the algebraic boundary of the amplituhedron. It consists of a union of Schubert divisors coming from the faces of a cyclic polytope. 
\end{exa}

We should mention that, in recent work \cite{brown2025positive}, Brown and Dupont formalize the construction of canonical differential forms using mixed Hodge theory. Let $X$ be an $n$-dimensional projective variety and fix a divisor $Y \subset X$ containing the singular points of $X$. We assume $(X,Y)$ to have \emph{genus zero} \cite[Section~3]{brown2025positive}. The divisor $Y$ plays the role of the hyperplane and Schubert arrangements in Example \ref{ex:introduction}. Brown and Dupont identify a canonical linear map sending a relative cycle $\sigma \in H_n(X, Y)$ (e.g., the polytope $P$) to a logarithmic $n$-form $\omega_\sigma \in \Omega^n_{\rm log}(X \setminus Y)$ (e.g., $\Omega(P)$). This suggests a new definition, and a new mechanism for constructing canonical forms. The ingredients remain roughly the same. 

Example \ref{ex:introduction} showcases a direct connection between positive geometry and (real) algebraic geometry via complex varieties, semi-algebraic sets and their algebraic boundary divisors. The varieties $X$ seen in physics are typically of special interest to algebraic geometers as well. For instance, $X$ is often a moduli space, such as ${\rm Gr}_{\mathbb{C}}(k,n)$ or $\overline{\cal M}_{0,n}$ \cite{Lam}. On the other hand, hyperplane arrangements, polytopes, and Schubert arrangements in the Grassmannian should appeal to combinatorialists. Indeed, studying interesting objects in positive geometry requires tools from several established areas of mathematics, including, but not limited to, algebraic geometry and combinatorics. Other relevant areas include tropical geometry, (algebraic) analysis, and numerical mathematics.

The name \emph{positive geometry} is slightly deceiving in that it suggest that one is only interested in real, positive points. In fact, a feature of the subject is that one  often combines the complex, real and tropical viewpoint on each of its objects. A key example is the Grassmannian ${\rm Gr}_{\mathbb{C}}(k,n)$. The positive part of ${\rm Gr}_{\mathbb{R}}(k,n)$ is an important instance of a positive geometry. It projects linearly to the
amplituhedron ${\cal A}$, and the images of its positroid cells form the analogue of a triangulation of ${\cal A}$. The {\em tropical Grassmannian} ${\rm Trop}({\rm Gr}_{\mathbb{C}}(2,n))$ is the space of trees with $n$ labeled leaves.
These trees correspond to tree-level Feynman diagrams in particle scattering. Tropical Grassmannians with $k > 2$ are only understood for small $n$.  Their relevance in physics led to {\em CEGM theory}  \cite{cachazo2019scattering}. 

\begin{figure}[h]
\centering
 \includegraphics[width = 7.5cm]{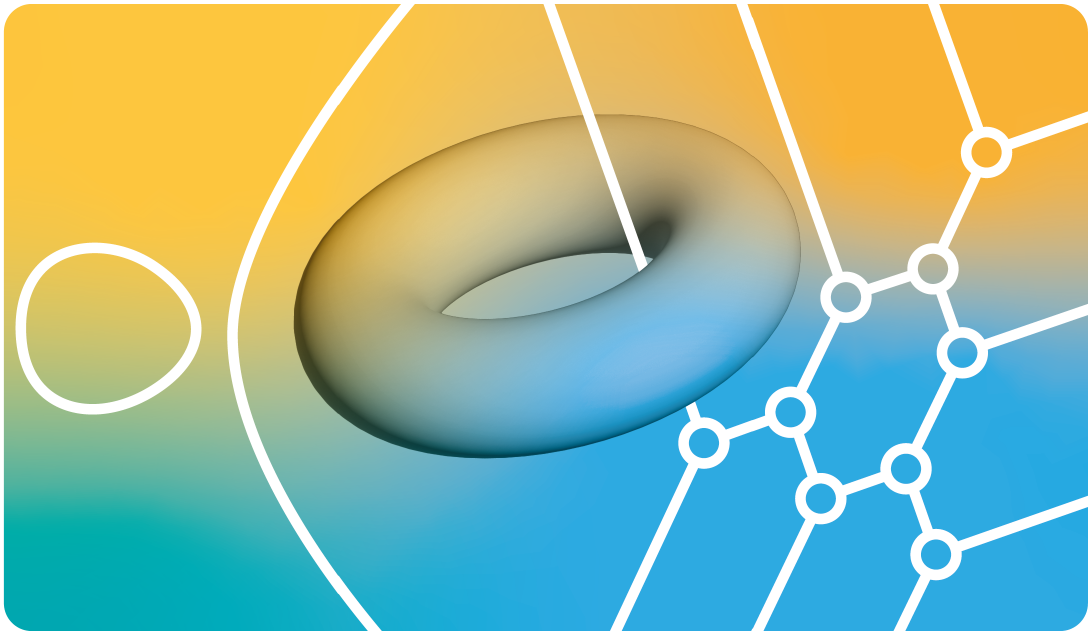}
\vspace{-0.02in}
\caption{Positive geometry rests on three flavors of
algebraic geometry: {\em real}, {\em complex} and {\em tropical}.
The diagram shows an elliptic curve in
 these three flavors.
 }
\label{fig:endler}
\end{figure}

A more visual example is the elliptic curve.
Its associated Riemann surface is a torus,
shown in the middle of Figure \ref{fig:endler}.
The cycle that indicates a genus one object
can be seen both in the left curve (real)
and in the right curve (tropical).

Though the inspiration for positive geometry comes from physics, its structures appear in other contexts too. For instance, a convex positive region $P$ represents the feasible set of an optimization problem. Here, the canonical dual volume function ${\rm vol}(P-x)^\circ$ from Example \ref{ex:introduction} is the universal barrier function of $P$. In algebraic statistics, the positive points $\mathbb{P}^n_{\geq 0}$ of projective space 
$\mathbb{P}^n$ are identified with the $n$-dimensional probability simplex; the space of probability distributions on a random variable with $n+1$ states. A statistical model is a complex projective variety $X \subset \mathbb{P}^n$, whose intersection with the probability simplex $X_{\geq 0} = X \cap \mathbb{P}^n_{\geq 0}$ is nonempty. Marginal likelihood is computed from Bayesian integrals on $X_{\geq 0}$ whose structure is similar to that of Feynman integrals. These settings justify the expectation that positive geometry will be widely applicable.

In summary, positive geometry revolves around  positive regions $P = X_{\geq 0}$ in complex algebraic varieties $X$, their boundary arrangements in $X$, existence and uniqueness of the canonical differential form $\Omega(X_{\geq 0})$, and the associated period integrals $\int_{\Gamma}\Omega(X)$. It uses a variety of techniques from combinatorics and complex, real and tropical geometry. We reiterate that this description and the articles in this volume form a snapshot of the field in its current nascent stage.

In what follows, we explain some of the concepts above in more detail. Our exposition is guided by examples and by the
contributions we are referring to.
We highlight connections between the cited papers, and
we explain how they fit the scope of positive geometry. The material is subdivided into three sections:
\begin{itemize}[leftmargin=0.6cm]
    \vspace{-0.1cm}\item Scattering equations and binary geometries \cite{Antolini,Betti,Bossinger,Cox,Kayser,Lam, Mazzucchelli, Rajan},
    \vspace{-0.2cm}\item (Semi-)algebraic subsets of the Grassmannian \cite{Agostini,Dian,ElMaazouz,Koefler,Mazzucchelli,Pratt,Seemann,De},
    \vspace{-0.2cm} \item Integrals and differential equations \cite{Borinsky,Fevola,Lam,Lotter,Pfister}.
\end{itemize}
Our article concludes with a collection of open problems and future directions.

\section{Scattering equations and binary geometries} \label{sec:2}

The framework of Cachazo--He--Yuan (CHY) for computing scattering amplitudes in bi-adjoint scalar $\phi^3$ theories illustrates many prominent concepts in positive geometry. The lecture notes \cite{Lam} provide an excellent mathematical exposition of this framework. It also contains open questions, some of which are answered in this volume \cite{Agostini,Bossinger}. This section presents a sketch of the CHY formalism. We work with $n = 5$ particles for concreteness, offering insight into what happens for general $n$. The $n = 5$ \emph{CHY amplitude} is the rational function
\begin{equation} \label{eq:ampln=5} A_5 \,\, = \,\,
 \frac{1}{s_{12}s_{45}} \, + \,  \frac{1}{s_{23}s_{15}} \, + \,  \frac{1}{s_{12}s_{34}} \, + \,  \frac{1}{s_{23}s_{45}} + \frac{1}{s_{15}s_{34}} \end{equation}
in the \emph{Mandelstam variables} $s_{ij} = s_{ji}$, $1 \leq i , j \leq 5$, with $s_{ii} = 0$. The $s_{ij}$ depend on the momenta of five massless particles. They satisfy $\sum_{j} s_{ij} = 0$ for all $i$ because of momentum conservation. Additionally, there can be a rank constraint on the 
 symmetric matrix $(s_{ij})_{i,j}$ with zeros on the diagonal. These conditions define the \emph{Mandelstam variety} of physically meaningful data $s_{ij}$. Different representations of these kinematic data lead to other determinantal varieties, studied in~\cite{Rajan}. 

Each term in \eqref{eq:ampln=5} corresponds to a triangulation of the pentagon or, dually, to a planar trivalent tree with five labeled leaves. The moduli space of planar trees is the positive tropical Grassmannian ${\rm Trop}^+  {\rm Gr}(2,n)$, and \eqref{eq:ampln=5} can be viewed as a sum over its maximal cones. Motivated by the analog of this observation for more general CEGM amplitudes, \cite{Antolini} studies the \emph{chirotropical Grassmannian}.

\begin{figure}
    \centering
    \includegraphics[width = 10cm]{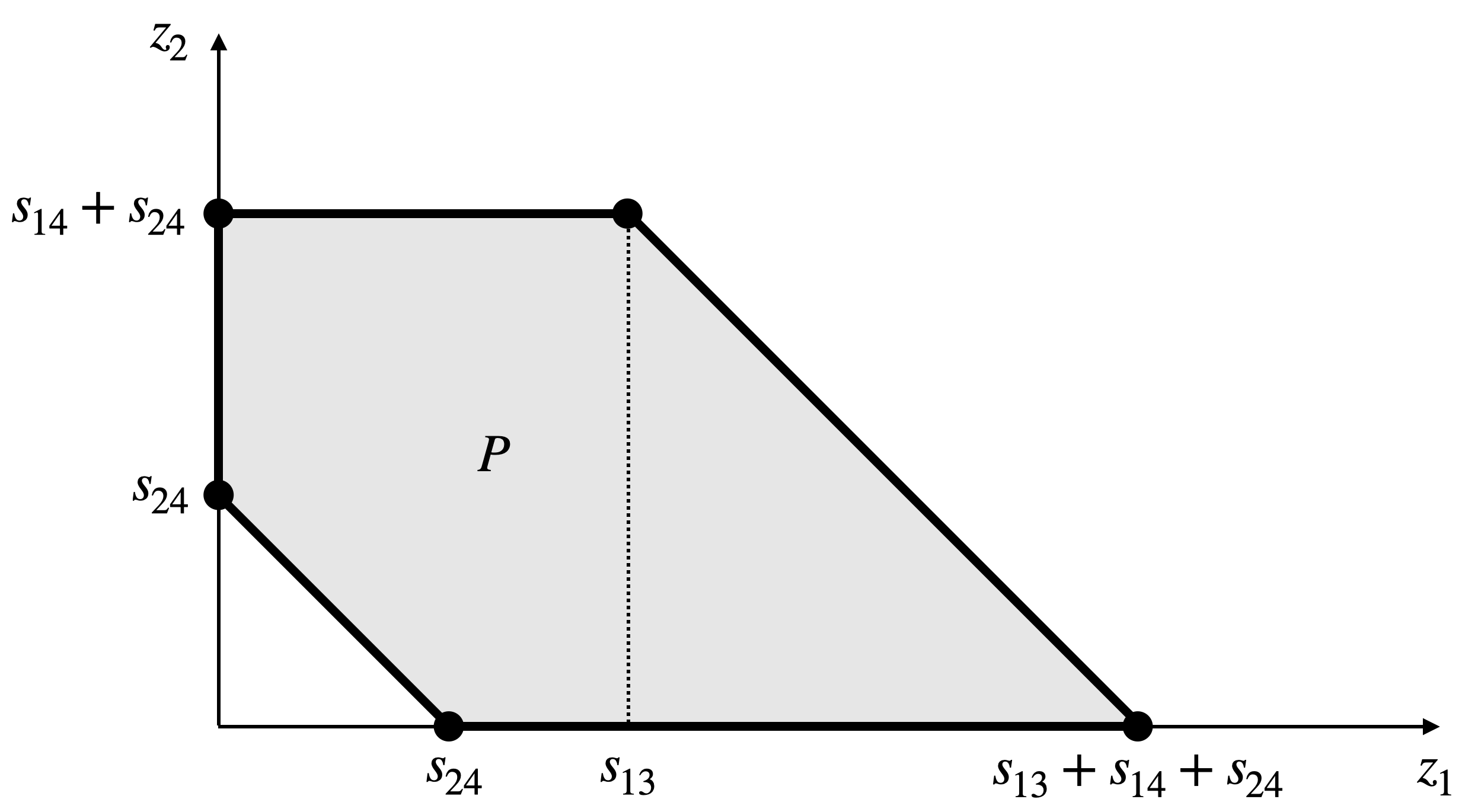}
    \caption{The amplitude $A_5$ is the dual volume of an associahedron.}
    \label{fig:ABHY}
\end{figure}

We now relate the function $A_5$ to the positive geometries from \cite{arkani2017positive} mentioned in the Introduction. We assume that $s_{13}$, $s_{14}$ and $s_{24}$ are positive real numbers. Let $X = \mathbb{P}^2$ and take $X_{\geq 0} = P$ to be the pentagon in Figure \ref{fig:ABHY}. 
Following Example \ref{ex:introduction}, its canonical form is $\Omega(P) = {\rm vol}(P-z)^\circ {\rm d} z$. We claim that substituting $z = (z_1,z_2) = (-s_{23},-s_{34})$ into the canonical function ${\rm vol}(P-z)^\circ$ gives \eqref{eq:ampln=5}. The pentagon $P$ is a realization of the two-dimensional associahedron \cite{arkani2018scattering}, whose vertices are the triangulations mentioned in the previous paragraph. 

The CHY formula \eqref{eq:CHYformula} writes \eqref{eq:ampln=5} as a sum over the complex critical points of a logarithmic potential function. Let $(p_{ij})_{i<j}$ be the $2 \times 2$ minors of 
\begin{equation}
\label{eq:twobyfive}
\begin{pmatrix}
1 & 1 & 1 & 1 & 0\\
0 & 1 & 1+x & 1+x+y & 1
\end{pmatrix}. 
\end{equation}
For $x,y \in \mathbb{C}^2 \setminus \{ \prod_{i<j}p_{ij} = 0 \}$, this matrix represents a point on the moduli space ${\cal M}_{0,5}$ of configurations of five distinct points $\{0,1,1+x,1+x+y,\infty\}$ on $\mathbb{P}^1$. That is, we model ${\cal M}_{0,5}$ as the complement of a \emph{hyperplane arrangement} $\{ (x,y) \, : \, \prod_{i<j}p_{ij} = 0 \} \subset \mathbb{C}^2$. The \emph{scattering equations} on ${\cal M}_{0,5}$ are given~by 
\begin{equation} \label{eq:scatteringeqs} \frac{\partial L}{\partial x} \, = \, \frac{\partial L}{\partial y} \, = \, 0 \quad \text{ where } \quad L \, = \,  \log \prod_{i<j} p_{ij}^{s_{ij}}. \end{equation}
The paper \cite{Betti} solves such equations by compactifying ${\cal M}_{0,5}$ and using numerical Gr\"obner degenerations. The \emph{logarithmic discriminant} is a polynomial in the $s_{ij}$ characterizing when two solutions collide. This is studied in \cite{Kayser}. In our example, it has a nice symmetric determinantal expression \cite[Example 8.1]{Kayser}.

Denote the Euler operators in $x$ and $y$ by $\theta_x = x \,  \frac{\partial}{\partial x}$ and $ \theta_y = y \, \frac{\partial}{\partial y}$. We have 
\begin{equation} \label{eq:CHYformula} A_5 \,\quad =  \sum_{(x,y) \in {\rm Crit}(L)}\!\! \det \begin{pmatrix}
    \theta_{x}^2 L & \theta_{x} \theta_{y} L \\ \theta_{x} \theta_{y} L & \theta_y^2 L
\end{pmatrix}^{\!\! -1}.
\end{equation}
The sum is over the set ${\rm Crit}(L)$ of solutions to \eqref{eq:scatteringeqs}. Hence, the coordinates $(x,y)$ are algebraic functions in $s_{ij}$. In our case, that set consists of only two points. The reader is encouraged,
by \cite[Exercise 3.40]{Lam}, to verify the identity \eqref{eq:CHYformula}.

We will reformulate \eqref{eq:scatteringeqs} in different coordinates on ${\cal M}_{0,5}$. The cross ratios
\[ [ij|kl] \, = \, \frac{p_{ik} p_{jl}}{p_{il}p_{jk}}, \quad \{i,j\}, \{k,l\} \in \binom{[5]}{2}, \, \,  \{i,j\} \cap \{k,l\}  = \emptyset \]
are regular functions on ${\cal M}_{0,5}$. The \emph{dihedral coordinates}, as introduced by Francis Brown \cite[Section 2]{Brown2009}, are the cross ratios $[i,i+1|j,j+1]$. These are labeled by diagonals $(ij)$ of a convex pentagon which are not edges. We denote these diagonals by $\Delta(5) = \{(13),(14),(24),(25),(35)\}$.  The map 
\[ {\cal M}_{0,5} \, \longrightarrow \, (\mathbb{C}^*)^{5}, \quad (x,y) \longmapsto ([i,i+1|j+1,j])_{(ij) \in \Delta(5)}.\]
is a closed embedding. The image $U \simeq {\cal M}_{0,5}$ is defined by five equations: 
\[ 
\begin{matrix} u_{13} + u_{24}u_{25} \, = \, 1, &
u_{24} + u_{13}u_{35} \, = \, 1, &
u_{35} + u_{14}u_{24} \, = \, 1, \\
u_{14} + u_{25}u_{35} \, = \, 1, &
u_{25} + u_{13}u_{14} \, = \, 1.
\end{matrix}
\]
Its closure in $\mathbb{C}^5$ is an instance of a \emph{binary geometry} \cite[Section 2.1]{Lam}. The paper \cite{Bossinger} discovers new binary geometries arising from \emph{pellytopes}. In \cite{Cox}, the authors study tropicalizations of binary geometries associated with cyclohedra. 

The potential function $L$ is given in dihedral coordinates by 
\[ L \, = \, \log \prod_{(ij) \in \Delta(5)} u_{ij}^{X_{ij}}, \quad \text{where} \quad X_{ij} \, = \, s_{i,j+1} + s_{i+1,j} - s_{ij} - s_{i+1,j+1}. \]
One checks that its critical points on $U$ satisfy the \emph{scattering equations}
\[ X^\top \footnotesize \begin{pmatrix}
0 & 1-u_{14}-u_{25} & 1-u_{13} & u_{13}-1 & u_{24} + u_{35}-1 \\
u_{14} + u_{35}-1 & 0 & 1-u_{13}-u_{25} & 1-u_{24} & u_{24}-1 \\
u_{35}-1 & u_{14} + u_{25} - 1 & 0 & 1-u_{13}-u_{24} &1-u_{35} \\
1-u_{14} & u_{14}-1 & u_{13}+u_{25}-1 & 0 & 1-u_{24}-u_{35} \\ 
1-u_{14}-u_{35} & 1-u_{25} & u_{25}-1 & u_{13}+u_{24}-1 & 0
\end{pmatrix}  \, = \, 0, \]
where $X^T$ is the row vector $(X_{13},X_{24},X_{35},X_{14},X_{25})$. It would be interesting to 
study scattering equations in dihedral coordinates for general $n \geq 5$. 

\section{(Semi-)algebraic subsets of the Grassmannian} \label{sec:3}

As mentioned in the Introduction, the discovery of the \emph{amplituhedron} in \cite{arkani2014amplituhedron} lies at the origin of positive geometry. Its canonical differential form is the integrand of the scattering amplitude in $N = 4$ Super Yang-Mills theory. The contribution of \cite{De} is a detailed explanation of this 
assertion, aimed at mathematicians.

Fix integer parameters $k, m, n$ satisfying $0 < k \leq k + m \leq n$,
and let $Z \in \mathbb{R}^{n \times (k + m)}$ be a real matrix, all of whose ordered $(k + m)\times(k+m)$ minors are positive. The map $\tilde{Z}: {\rm Gr}(k,n)_{\geq 0} \rightarrow {\rm Gr}(k,k+m)$ is that induced by matrix multiplication: $[X] \mapsto [X \cdot Z]$ for $X \in \mathbb{R}^{k \times n}$ of rank $k$ with non-negative $k \times k$ minors.

The \emph{tree-level amplituhedron} ${\cal A}_{n,k,m}(Z)$ is the image of the map $\tilde{Z}$: 
\[ {\cal A}_{n,k,m}(Z) \, = \, \tilde{Z}({\rm Gr}(k,n)_{\geq 0}) \, \subset \,  {\rm Gr}(k,k+m).\]
It is a semi-algebraic subset of the real locus ${\rm Gr}_{\mathbb{R}}(k,k+m) $ in $ {\rm Gr}(k,k+m)$.  

\begin{figure}
    \centering
    \includegraphics[width = 7cm]{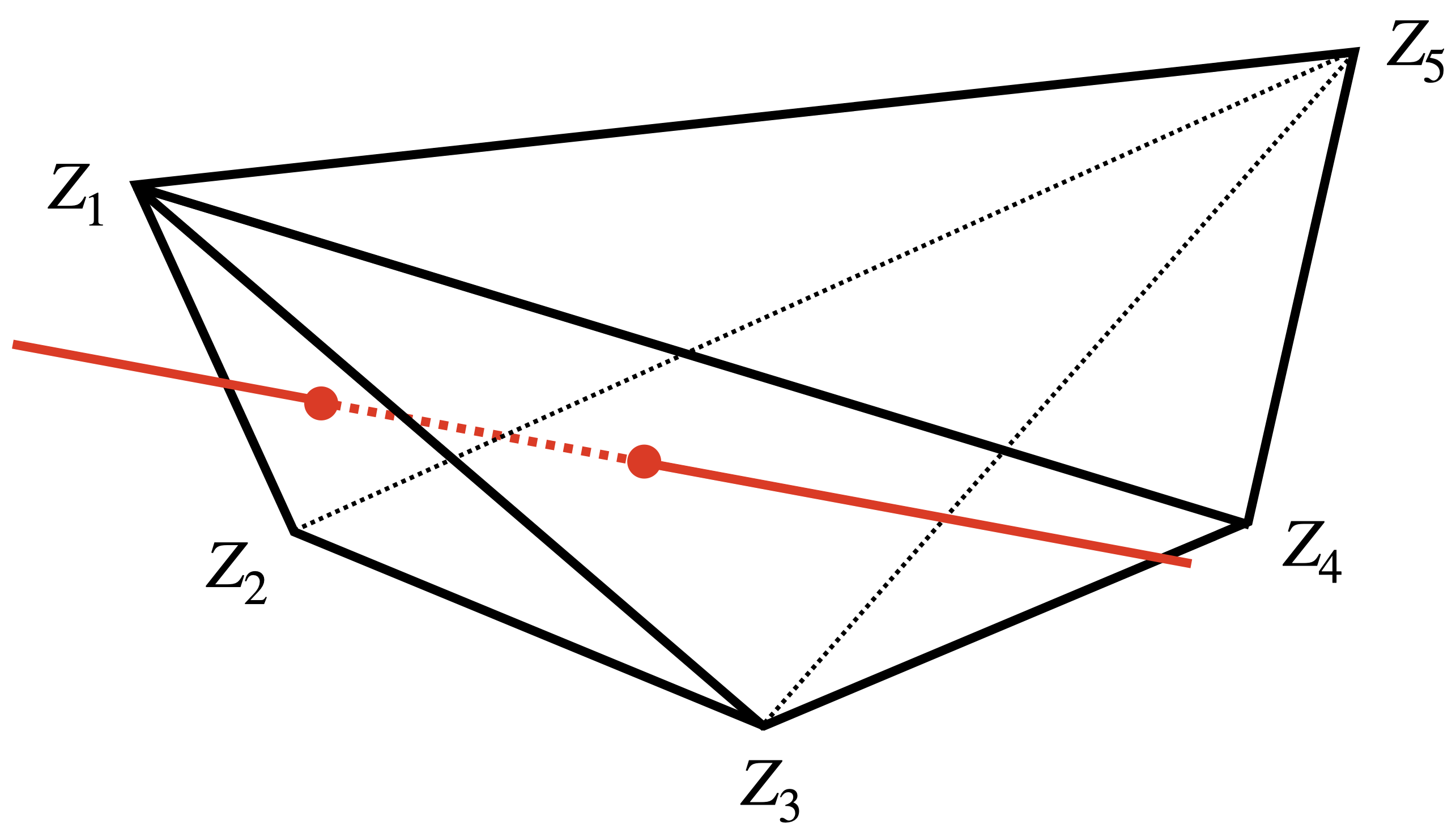}
    \caption{A cyclic polytope in $\mathbb{RP}^3$ with a stabbing line.}
    \label{fig:cyclic}
\end{figure}

For a concrete example, set $k = m = 2$, $n = 5$. To make sure all $4 \times 4$ minors of $Z \in \mathbb{R}^{5 \times 4}$ are positive, we can pick points on the twisted cubic for its rows: 
\begin{equation} \label{eq:myZ} Z_i \,=\, \begin{pmatrix}  1 & i & i^2 & i^3
\end{pmatrix}, \quad i = 1, \ldots, 5. \end{equation}
These are the vertices of a cyclic polytope ${\cal A}_{5,1,3}(Z) = {\rm Conv}(Z)$ in $\mathbb{RP}^3$.
The combinatorial type of this polytope is shown in Figure \ref{fig:cyclic}. 
Since ${\cal A}_{5,2,2}(Z)$ is a full-dimensional semi-algebraic subset of ${\rm Gr}_{\mathbb{R}}(2,4)$, it is described by polynomial inequalities in the Pl\"ucker coordinates $p_{12}, p_{13}, \ldots, p_{34}$ on ${\rm Gr}(2,4)$. This description uses the standard notation $\langle A B i j \rangle$ for the linear form in $p_{ij}$ vanishing on lines $AB \in {\rm Gr}(2,4)$ which intersect the line $Z_iZ_j$. It was conjectured in \cite{arkani2018unwinding}, and proved in \cite{parisi2023m}, that the lines $AB \in {\cal A}_{5,2,2}(Z)$ are characterized as follows:
\[ \begin{matrix}
    \langle AB12 \rangle , \, \langle AB23 \rangle  ,  \langle AB34 \rangle  , \, \langle AB 45 \rangle, \, \langle AB15 \rangle \, \, \text{ have the same sign,}\\
    \text{the sequence }\langle AB12 \rangle, \, \langle AB13 \rangle, \, \langle AB14 \rangle, \, \langle AB15 \rangle \, \text{ has exactly two sign flips}.
\end{matrix}\]
Here zeros are ignored when counting sign flips. The red line $AB$ in Figure \ref{fig:cyclic} passes through the interior of the facets $Z_1Z_2Z_3$ and $Z_1Z_3Z_4$. It belongs to the \emph{$2$-stabbing set} of ${\rm Conv}(Z)$, but it does not belong to the amplituhedron: $\langle AB12 \rangle$ and $\langle AB34 \rangle$ have opposite signs. However, one can show that each line in the amplituhedron stabs the cyclic polytope. The paper \cite{Seemann} introduces and studies more general stabbing sets and their relation to tree amplituhedra. 

The algebraic boundary of ${\cal A}_{5,2,2}(Z)$ is a union of five hyperplane sections of ${\rm Gr}(2,4)$, given by the equation $\langle AB12 \rangle \langle AB23 \rangle \cdots \langle AB15 \rangle = 0$. Such arrangements in the Grassmannian are natural generalizations of hyperplane arrangements in $\mathbb{P}^n$, which are extensively studied in (algebraic) combinatorics and which we encountered in Section \ref{sec:2}. The paper \cite{Mazzucchelli} investigates the real and complex topology of hyperplane arrangement complements in ${\rm Gr}(k,n)$. 

The \emph{adjoint hypersurface} of our tree-level amplituhedron ${\cal A}_{5,2,2}(Z)$, 
for the specific $Z$ in \eqref{eq:myZ}, is given by a linear equation in Pl\"ucker coordinates \cite{Ranestad_Sinn_Telen_2024}: 
\begin{equation} \label{eq:adjA5} 593 \, p_{12} - 330 \, p_{13} + 143 \, p_{14}  + 49 \, p_{23} - 30 \, p_{24} + 5 \, p_{34} \, = \, 0. \end{equation}
It is uniquely determined by the following interpolation conditions. The adjoint contains the unique line $(1|23|45) \in {\rm Gr}(2,4)$ passing through $Z_1$ and intersecting the lines $Z_2Z_3$ and $Z_4Z_5$. Similarly, it contains $(2|34|15)$, $(3|45|12)$, $(4|15|23)$ and $(5|12|34)$. The linear form \eqref{eq:adjA5} is the numerator of the integrand of the corresponding $N=4$ Super Yang-Mills amplitude.
The authors of \cite{Koefler} seek to uniquely determine the numerator of \emph{MHV gravity amplitudes} via similar interpolation conditions on a flag variety. Except for some small cases, this description remains conjectural. It is natural to ask whether this numerator is adjoint to a positive geometry -- an ``amplituhedron'' for MHV amplitudes.

Tree amplituhedra are a family of candidate positive geometries inside $X = {\rm Gr}(k,k+m)$ by considering $X_{\geq 0} = {\cal A}_{n,k,m}(Z), \, n \geq k + m$. When $k = 1$, ${\cal A}_{n,k,m}(Z)$ is a cyclic polytope in $\mathbb{P}^{k + m - 1}$, and hence a positive geometry (Example \ref{ex:introduction}). When $k + m = n$, ${\cal A}_{n,k,m}(Z)$ is isomorphic to ${\rm Gr}(k,n)_{\geq 0}$, which is known to be a positive geometry. Recent progress for $k = m = 2$ is found in \cite{Ranestad_Sinn_Telen_2024}. For other $n,k,m$, the amplituhedron ${\cal A}_{n,k,m}(Z)$ is only conjectured to be a positive geometry \cite[Conjecture 9]{arkani2017positive}. \emph{Loop-level amplituhedra}, as opposed to tree-level amplituhedra, are known to \emph{not} be positive geometries in the sense of \cite{arkani2017positive}. They still have a meaningful adjoint hypersurface, as illustrated by the case study~\cite{Dian}. 

Other interesting positive geometries are semi-algebraic subsets of subvarieties of the Grassmannian. This arises in \emph{ABJM scattering amplitudes}, with the relevant subvariety being the \emph{orthogonal Grassmannian}. The authors of \cite{ElMaazouz} examine the \emph{positive orthogonal Grassmannian}, and prove that it is a positive geometry in some cases. A different subvariety, called \emph{ABCT variety}, is studied in \cite{Agostini}, which solves part of \cite[Problem 4.18]{Lam}. Lam's conjecture \cite[Conjecture 4.10]{Lam} regarding a positive geometry inside the ABCT variety remains~open. 

Inspired by the amplituhedron map $\tilde{Z}$, the paper \cite{Pratt} investigates whether a $(km-1)$-dimensional subvariety of ${\rm Gr}(k,n)_{\geq 0}$ can be recovered from its linear projections $[X] \mapsto [X \cdot Z] \in {\rm Gr}(k,k+m)$ for all possible $Z$. This generalizes the classical construction of recovering a projective variety from its \emph{Chow form}.

\section{Integrals and differential equations} \label{sec:4}

Observables in physics are usually expressed as integrals.
Such integrals include  scattering amplitudes in particle physics
and correlators in cosmology.  The following Euler integral
serves as a blueprint for cosmological correlators:
$$ \phi(c) \,\,=\,\, \int_\Gamma  \frac{1}{ (c_1 \alpha_1 + c_2 \alpha_2 + c_3)
(c_4 \alpha_1 + c_5) (c_6 \alpha_2 + c_7)} 
\alpha_1^{\epsilon+1} \alpha_2^{\epsilon+1}  
\frac{{\rm d} \alpha_1}{\alpha_1} \frac{{\rm d} \alpha_2}{\alpha_2} .
$$
Here $\Gamma$ is a twisted cycle in $\mathbb{C}^2$.
The integral is a function $\phi$ of
$ c = (c_1,c_2, \ldots , c_7 )$. It satisfies
a system of linear differential equations known as the
{\em GKZ system}. These equations come in two flavors.
First, $\phi$ is annihilated by Euler operators
\begin{equation}
\label{eq:euler} \begin{matrix} c_1 \partial_1 + c_2 \partial_2 + c_3 \partial_3 + 1\, ,\,\,
c_4 \partial_4 + c_5 \partial_5 + 1\,,\,\,
c_6 \partial_6 + c_7 \partial_7 + 1\,, \\
c_1 \partial_1 + c_4 \partial_4 + (\epsilon+1) \,,\quad
c_2 \partial_2 + c_6 \partial_6 + (\epsilon+1).
\end{matrix}
\end{equation}
Here, $\partial_i $ denotes the operator $\frac{\partial}{\partial c_i}$.
And, second, we have the toric operators
\begin{equation}
\label{eq:toric}
\partial_1 \partial_5 - \partial_3 \partial_4 \quad {\rm and} \quad
\partial_2 \partial_7 - \partial_3 \partial_6.
\end{equation}
The operators in (\ref{eq:euler})-(\ref{eq:toric}) 
generate a holonomic $D$-ideal of rank $4$.
This means in particular that the solution space is $4$-dimensional.
This GKZ system and its role for the
two-site chain model in cosmology are explored in
\cite[Section~3.1]{Fevola}. 

In most applications, the coefficients $c_i$ of the parenthesized linear functions
are not independent unknowns, but they satisfy special constraints.
This makes the derivation of annihilating differential operators more
difficult. One scenario, of relevance to cosmology, is studied
in the article \cite{Pfister}. Here, the linear forms are fixed,
and only the constant terms are allowed to vary. For instance, consider
$$ \psi(u) \,\,=\,\, \int_\Gamma  ( \alpha_1 +  \alpha_2 + u_1)^{s_1}
( \alpha_1 + u_1)^{s_2} ( \alpha_2 + u_3)^{s_3}
\alpha_1^{\nu_1} \alpha_2^{\nu_2}
\frac{{\rm d} \alpha_1}{\alpha_1} \, \frac{{\rm d} \alpha_2}{\alpha_2} . $$
The five exponents $s_1,s_2,s_3,\nu_1,\nu_2$ are parameters.
The $D$-module for the function $\psi(u)$ in three
variables $u = (u_1,u_2,u_3)$ is worked out
 in \cite[Section 4.5]{Pfister}.

The CHY amplitude we saw in (\ref{eq:ampln=5}) also has 
an analytic representation by means  of an integral.
 Namely, $A_5$ arises as a limit 
from the {\em string amplitude} 
\begin{equation*}  \phi_\epsilon(s) \,\,=\,\, \epsilon^2 \int_{\mathbb{R}^2_+} \exp(\epsilon \cdot L) \, \frac{{\rm d}x}{x} \, \, \frac{{\rm d}y}{y} \!\, = \, \,\epsilon^2 \int_{\mathcal{M}_{0,5}^+} \frac{1}{p_{12} p_{23} p_{34} p_{45} p_{51}}
\prod_{1 \leq i < j < 5} \!p_{ij}^{\epsilon \cdot s_{ij}} \,{\rm d}p. \end{equation*}
Here $L(x,y)$ is the scattering potential from \eqref{eq:scatteringeqs}, which depends on the Mandelstam variables $s = (s_{12}, \ldots, s_{45})$. The positive moduli space $\mathcal{M}_{0,5}^+$ is the region $\mathbb{R}^2_{+}$ in the hyperplane arrangement mentioned prior to (\ref{eq:scatteringeqs}).
The integrand is the {\em Koba--Nielsen potential} from \cite[Section 5.1]{Lam}.
The limit of $\phi_\epsilon(s)$, as $\epsilon$ tends to zero,
is the rational function $A_5$ in \eqref{eq:ampln=5} or 
\eqref{eq:CHYformula}. This is the content of \cite[Theorem 5.12]{Lam}.
In string theory, the parameter $\epsilon = \alpha'$ is the \emph{inverse string tension}.

{\em Iterated integrals} arise in particle physics, number theory, 
and many other settings. In stochastic analysis, these are defined for paths
$\, X: [0,1] \rightarrow \mathbb{R}^d$. Every sequence
$i_1,i_2,\ldots,i_k$ of indices in $\{1,2,\ldots,d\}$ specifies an iterated integral
$$ \sigma_{i_1 i_2 \cdots i_k}(X) \,\,\, = \,\, \int_{\Delta_k} {\rm d} X_{i_1}(t_1)\, {\rm d}X_{i_2}(t_2) \,\cdots \, {\rm d}X_{i_k}(t_k). $$
Here  $\Delta_k$ denotes the $k$-dimensional simplex given
 by $0 \leq t_1 \leq t_2 \leq \cdots \leq t_k \leq 1$.
These real numbers are the entries in a $k$-dimensional tensor
$\sigma$ of format $d \times d \times \cdots \times d$.
This tensor is {\em $k$-th signature} of the path $X$.
The article \cite{Lotter} concerns piecewise linear paths $X$, whose
breakpoints are the vertices of a cyclic polytope. It characterizes
linear combinations of the tensor entries $\sigma_{i_1 i_2 \cdots i_k}(X)$  
that express invariant,
geometric properties of the polytope.
Cyclic polytopes are important in particle physics
because they are special cases of amplituhedra.

The article \cite{Borinsky} grew out of the lecture series
{\em Counting graphs using quantum field theory},
given by Michael Borinsky
at the May 2024 school in Leipzig.
The  integrals studied in that article are similar to those
seen in quantum field theory:
$$ I(z) \,\,\, = \,\,\,\frac{z}{2 \pi} \int_D {\rm exp}(z \,g(x,y)) \,{\rm d} x \,{\rm d} y. $$
Here, $D$ is a subset of $\mathbb{R}^2$
and $g : D \rightarrow \mathbb{R}$ is an appropriate function. Notice the resemblance to the string amplitude $\phi_\varepsilon(s)$, setting $g = L$ and $z = \epsilon$.
The asymptotic expansion of $I(z)$ for large $z$ serves
as a generating function for the number of regular edge-bicolored graphs.
Methods from analysis and combinatorics are combined to determine
the limit behavior of these numbers. Analogous results 
in one variable had huge implications for the cohomology of moduli spaces,
both classical and tropical, and we are similarly optimistic for the impact of \cite{Borinsky}.

 In the Introduction,   we emphasized
 the physical importance of the volume function $ \,   x \mapsto {\rm vol}(P-x)^\circ $.
 As in Example \ref{ex:introduction},
  $P$ is a convex polytope or a more general positive geometry.
 Of course, such a volume is an integral. If $P$ depends on
 parameters, then so does the function $x \mapsto {\rm vol}(P-x)^\circ$.
 It therefore makes sense  to apply the {\em differential equations method},
 which represents such a  holonomic function by a  $D$-module.
The study of integrals and differential equations is thus an essential
ingredient for the emerging field of positive geometry.

\section{Open questions}

In this section we present one open problem
 for each of the articles in this collection.
In each case, 
the problem we chose is either stated in
the article, in which case a  direct citation is given, or it is inspired by results from the article.

\medskip \noindent {\bf Problem} \cite{Agostini}: 
The image of the positive Grassmannian  ${\rm Gr}(2,n)_{\geq 0}$
in the ABCT variety  $V(d+1,n)$ is a semialgebraic set. Is that image a positive geometry?
Is $V(d+1,n)$ with the relevant boundary a genus zero pair, in the sense of
\cite{brown2025positive}?

\medskip \noindent {\bf Problem} \cite{Antolini}:
Consider a uniform oriented matroid $M$ whose realization space 
is diffeomorphic to an open ball. Is the chirotropical
Dressian of $M$  equal to the chirotropical 
Grassmannian of $M$?   \ 
This is conjectured in \cite[Conjecture~6.2]{Antolini}.

\medskip \noindent {\bf Problem} \cite{Betti}:
A weak formulation of the CHY scattering equations leads to the intersection of the reciprocal variety of the moduli space $\mathcal{M}_{0,n}$ with a linear space. 
Determine the scheme structure of this intersection \cite[Conjecture~5.6]{Betti}.

\medskip \noindent {\bf Problem} \cite{Borinsky}:
How to best generalize the results from this paper to
graphs with three or more colors?
Now we must integrate over a space $D$ of dimension $\geq 3$.

\medskip \noindent {\bf Problem} \cite{Bossinger}: 
Does the Pellspace lead to an interesting compactification of
the moduli space $\mathcal{M}_{0,n}$?
Develop this from the perspective of
algebraic geometry.

\medskip \noindent {\bf Problem} \cite{Cox}:
Does the space of axially symmetric phylogenetic trees coincide with the tropicalized cluster configuration space in type C?
\cite[Conjecture~7.5]{Cox}.

\medskip \noindent {\bf Problem} \cite{De}:
Is the amplituhedron a positive geometry, in the sense of \cite{arkani2017positive}?
This is a guiding conjecture, going back to the original paper by Arkani-Hamed and Trnka 
\cite{arkani2014amplituhedron}. See \cite[Section 6]{De}.
Can the paper \cite{brown2025positive} shed new light on this?

\medskip \noindent {\bf Problem} \cite{Dian}:
Show that the two-loop amplituhedron $\mathcal{A}^{(2)}_4$ is a
weighted positive geometry. Following \cite[Section 5.2]{Dian},
some residues need to be computed.

\medskip \noindent {\bf Problem} \cite{ElMaazouz}:
Determine a  cell decomposition for the 
positive orthogonal Grassmannian. 
What should be added to the orthopositroids? \qquad
\cite[Problem 5.4]{ElMaazouz}.

\medskip \noindent {\bf Problem} \cite{Fevola}:
Prove the partial fraction decomposition 
suggested in \cite[Conjecture 3.5]{Fevola}
for the flat space wave function
$\psi_{\rm flat}$ of an arbitrary Feynman graph.

\medskip \noindent {\bf Problem} \cite{Kayser}:
Show that the logarithmic discriminant of the moduli space
$\mathcal{M}_{0,n}$ is an irreducible hypersurface, and
determine its degree. \quad
\cite[Conjecture 1]{Kayser}.

\medskip \noindent {\bf Problem} \cite{Koefler}:
Establish the unique representation of
the $n$-point MHV amplitude. This is stated
as a commutative algebra problem
in \cite[Conjecture 2.1]{Koefler}.

\medskip \noindent {\bf Problem} \cite{Lam}:
Study spinor-helicity polytopes, with focus on their
parametrizations and associated cohomology classes. Solve
\cite[Problems 4.22 and 4.23]{Lam}.

\medskip \noindent {\bf Problem} \cite{Lotter}:
Determine the ring of volume invariants. This refers to
signatures of piecewise-linear paths on cyclic polytopes.
Prove \cite[Conjecture~4.1]{Lotter}.

\medskip \noindent {\bf Problem} \cite{Mazzucchelli}:
What is the maximal number of regions in the complement
of an arrangement of $d$ Schubert hyperplanes 
in the real Grassmannian
${\rm Gr}_\mathbb{R}(k,n)$?

\medskip \noindent {\bf Problem} \cite{Pfister}:
Determine a set of generators for the
annihilator ${\rm Ann}_{D(s,V)}(\phi)$ of the correlator function $\phi$
of an arbitrary hyperplane arrangement $V$.

\medskip \noindent {\bf Problem} \cite{Pratt}:
Let $\mathcal{V}$ be a subvariety of the Grassmannian
${\rm Gr}_\mathbb{C}(k,n)$ which is
defined by $k(n-r)+1$ generic linear forms in the Pl\"ucker coordinates.
Prove that $\mathcal{V}$ coincides with its
Chow-Lam recovery $\mathcal{W}_\mathcal{V}$. \hfill
\cite[Conjecture 6.6]{Pratt}.

\medskip \noindent {\bf Problem} \cite{Rajan}:
Show that the kinematic variety for spinor brackets of order $\leq 3$
is the variety of tensors with multilinear rank $\leq (2,4,2)$. \ \ \
\cite[Conjecture 6.2]{Rajan}.

\medskip \noindent {\bf Problem} \cite{Seemann}:
Show that the cells in the $k$-face Schubert arrangement of
any convex polytope are determined by their sign vectors. \hfill
\cite[Conjecture 3.5]{Seemann}.

\medskip

We hope that these open problems will
inspire mathematical progress in the development of positive geometry,
and that this will have an impact in physics.

\bibliography{references}

\end{document}